# Multiple Normalized Solutions to a Class of Modified Quasilinear Schrödinger Equations

Zhouxin Li [a], Ayesha Baig*, [a]

[a] Department of Mathematics and Statistics, Central South University, Changsha 410083, PR China

**Abstract**: In our research, we focus on the existence, non-existence, and multiplicity of positive solutions to a Quasilinear Schrödinger equation in the form:

$$-\Delta u + \lambda u + \frac{k}{2}[\Delta(u^2)]u = f(u), \quad u \in H^1(\mathbb{R}^N)$$

With prescribed mass:

$$\int_{\mathbb{R}^N} |u|^2 dx = c,$$

Here $N \geq 3$, The dual approach is used to transform this equation into a corresponding semi-linear form. Then, we implement a global branch approach, adeptly handling nonlinearities $f(s)$ that fall into mass subcritical, critical, or supercritical categories. Key aspects of this study include examining the positive solutions' asymptotic behaviors as $\lambda \to 0^+$ or $\lambda \to +\infty$ and identifying a continuum of unbounded solutions in $(0, +\infty) \times H^1(\mathbb{R}^N)$.

1. Introduction

1.1. Background and Motivation

In this work, we engage in a comprehensive analysis of quasilinear Schrödinger equation (QLSE), a subject that has gained considerable traction in mathematical physics due to its complex nature and broad applications. Central to our study are the solutions which conform to the equation:

$$\begin{cases} -\Delta u + \lambda u + \frac{k}{2}[\Delta(u^2)]u = f(u), \\ \int_{\mathbb{R}^N} |u|^2 dx = c, \end{cases} \quad x \in \mathbb{R}^N \qquad (1.1)$$

where $\lambda$ appears as a Lagrange multiplier. This exploration, situated within the functional space:

___________________
*Corrosponding Author

$$H^1(\mathbb{R}^N) = \{u \in L^2(R^N): \nabla u \in L^2(\mathbb{R}^N)\}.$$

under the constraint:
$$Y(c) = \{u \in H^1(\mathbb{R}^N) \mid \int_{R^N} k|\nabla u|^2 u^2 \, dx < \infty, \|u\|_{L^2}^2 = c\}.$$

The exploration of solutions to (1.1) is intricately connected with the identification of standing wave solutions, which are essential for a comprehensive understanding of the behavior of the time-dependent quasilinear Schrödinger equation. This specific equation can be expressed as:

$$i\,\partial_t z = -\Delta z - l(|z|^2)z + \frac{\kappa}{2}[\Delta\rho(|z|^2)]\rho'(|z|^2)z, \quad (t,x) \in \mathbb{R}^+ \times \mathbb{R}^N \tag{1.2}$$

Here, $i$ symbolizes the imaginary unit. The function $z: \mathbb{R} \times \mathbb{R}^N \to \mathbb{C}$, $l$ and $\rho$ are real-valued functions.

Equations of this form, as described in (1.2), frequently emerge in the field of mathematical physics and have been applied as mathematical models in various physical contexts, particularly where the nonlinear term $\rho$ is involved. For instance, the specific case of $\rho(s) = s$ has been utilized in plasma physics to model the superfluid film equation, as discussed by Kurihara [6]. For $l(s) = as^p$, (1.2) appears in various problems in plasma physics and nonlinear optics, e.g. oscillating soliton instabilities during microwave and laser heating of plasma [36,37]. (1.2) is also the basic equation describing oscillations in a superfluid film when $(s) = -\alpha - \frac{\beta}{(a+s)^3}$. Normalized solutions are pivotal in modeling and understanding complex behaviors in various physical systems, such as optical fibers and quantum fluids, where nonlinear interactions are key. They offer critical insights into the dynamic stability and wave propagation characteristics in these systems, making them indispensable in both theoretical physics and practical applications but little is known about the existence (or non-existence) of normalized solutions compared with the problem where $\lambda$ is prescribed.

For the first time, mass super-critical case was considered in [5], since the energy functional is unbounded below on:
$$S_c = \{u \in H^1(R^N), \| u \|_2^2 = c\},$$

Subject to $S_c$ and $f(s)$ is a nonlinearity satisfy the following global condition.

(H1) $f: \mathbb{R} \to \mathbb{R}$ is continuous and odd.

(H2) $\exists (\alpha, \beta) \in \mathbb{R} \times \mathbb{R}$ satisfying
$$2 + \frac{4}{N} < \alpha \leq \beta < \frac{2N}{N-2},$$

Such that,
$$\alpha F(s) \leq f(s)s \leq \beta F(s),$$

With $F(s) = \int_0^s f(t)dt$.

(H3) $\tilde{F}: \mathbb{R} \to \mathbb{R}, \tilde{F} = f(s)s - 2F(s)$ such that:
$$\tilde{F}'(s)s > \frac{2N+4}{N}\tilde{F}(s).$$

Jeanjean's approach, using a minimax argument and smart compactness argument, proved instrumental in obtaining normalized solutions for these types of equations.

Building on this foundation, Bartsch and de Valeriola [15] expanded the understanding of these solutions' multiplicity, while Ikoma and Tanaka [16] provided alternative proofs using ideas related to symmetric mountain pass theorems. Later on [17] L. Jeanjen with Sheng-Sen Lu has revisted [5] and they were able to study the problem with less strict conditions on the nonlinearity.

More recently, research has ventured into more general nonlinearities, as explored in studies [18].

All this work is done on the semi-linear problem but for the quasilinear problem, Due to the quasilinearity the variational formulation is not smooth in general. The difficulties also arise for similar problems without constraints. Similar problems without constrains have been studied extensively in the last two decades (e.g. [23-27] and references therein) for which one has the following quasi-linear elliptic equation with $\mu$ being a fixed negative constant:
$$-\Delta u - \Delta(|u|^2)u - \mu u = |u|^{p-2}u, \quad \text{in } \mathbb{R}^N, \tag{1.3}$$

With the corresponding energy functional:
$$\bar{J}(u) = \frac{1}{2}\int_{\mathbb{R}^N} |\nabla u|^2 dx + \int_{\mathbb{R}^N} |\nabla u|^2 u^2 - \frac{\mu}{2}|u|^2 dx - \int_{\mathbb{R}^N} f(u)\,dx, \tag{1.4}$$

Here, $f(u) = |u(x)|^{p-2}u$. To overcome the non-differentiability of energy functional, dual approach was used and then minimax principle to obtain the existence of normalized solution. An important contribution in this area was made by an author [19], who considered (1.3) with a nonlinearity of the form $f(u) = |u(x)|^{p-2}u$ and used a perturbation method to obtain a ground state and an infinite number of normalized solutions. Tang and Yu's research [31] revisits the Berestycki-Lions conditions on ground state solutions for a nonlinear Schrödinger equation with variable potentials. Their innovative approach, leveraging the Jeanjean-Toland monotonicity trick and the IIP inequality, provides a fresh perspective on obtaining nontrivial critical points and least energy solutions, even when faced with non-sign definite functions.

One of the pivotal studies by Wei and Wu [28] addresses normalized solutions for Schrödinger equations with critical Sobolev exponent and mixed nonlinearities. Their work provides crucial insights into the compactness of minimizing sequences, leveraging the Aubin-Talenti bubbles and ground state solutions to construct test functions, a strategy hinting at the profound complexity of critical and supercritical perturbations. Quasilinear Schrodinger Equation with saturable Nonlinearity is studied in [35].

In 2023, Chergui, Gou, and Hajaiej [29] delve into the existence and dynamics of normalized solutions to nonlinear Schrödinger equations with mixed fractional Laplacians. Their work highlights the varied landscape of normalized solutions, from global minimizers in mass subcritical cases to saddle type critical points or local minimizers in mass critical or supercritical scenarios. The study underscores the shift in methodology required when transitioning from subcritical to critical or supercritical cases, marking a departure from global minimization problems to more complex minimax arguments. In the same year, A study presented in [30] on normalized solutions to Schrödinger equations with potential and inhomogeneous nonlinearities in large convex domains explores the existence and multiplicity of normalized solutions, adapting robust methods that can be applied to a range of other elliptic partial differential equations with potentials. This work, while focusing on a specific setting, broadens the application horizon for techniques used in analyzing normalized solutions.

Authors in [1] studied existence and nonexistence of nonzero solution of the following equation:

$$-\Delta u + V(x)u + \lambda u(x) + \frac{k}{2}[\Delta(u^2)]u = f(u), \quad x \in R^N,$$

with $V(x)$ as continuous potential by applying the dual approach and pohozaev identity with superlinear or asymptotically linear terms. This Equation is also considered in [32] with critical or supercritical exponents to study the existence of nontrivial solution. Here the author consider the continuous potential. More recently, L. Jeanjean, J. Zhang and X. Zhong [33] considered:

$$\begin{cases} -\Delta u + \lambda u = a(u), \\ \int_{\mathbb{R}^N} |u|^2 dx = c, \end{cases}, \quad x \in \mathbb{R}^N \quad (1.5)$$

and introduced a new approach to deal with problem (1.5) in the mass-subcritical, mass critical and mass-supercritical case in one unified way. In particular, the assumption (H1), (H2) and (H3) of [5] can be replaced by less strict assumptions:

$(A_1) a \in C^1[0, +\infty), a(0) = 0$ and $a(s) > 0$ for $s > 0$.

$(A_2)$ There exist $\alpha, \beta, \mu_1, \mu_2 > 0$ satisfying:

$$2 < \alpha, \beta < 2^* := \frac{2N}{N-2}.$$

such that,

$$\lim_{t \to 0^+} \frac{a'(t)}{t^{\alpha-2}} = \mu_1(\alpha - 1) > 0 \text{ and } \lim_{t \to +\infty} \frac{a'(t)}{t^{\beta-2}} = \mu_2(\beta - 1) > 0.$$

($A_3$) There is no positive, radially decreasing classical solutions for $-\Delta u = a(u)$ in $\mathbb{R}^N$.

Motivated by [33] we will extend the results to Quasilinear Schrödinger Equation of form (1.1). The examination of standing wave solution of the form:

$$z(t,x) = e^{-i\lambda t}u(x),$$

we observe that $z(t,x)$ satisfy Eq.(1.2) iff the function $u(x)$ solves (1.6) with $f(u) = l(u^2)u$.

$$-\Delta u + \lambda u + \frac{k}{2}[\Delta(u^2)]u = f(u), \quad x \in \mathbb{R}^N \tag{1.6}$$

And the natural energy functional associated to the equation is:

$$\bar{J}(u) = \frac{1}{2}\int_{\mathbb{R}^N}(1 - \kappa u^2)|\nabla u|^2 dx + \frac{1}{2}\int_{\mathbb{R}^N}\lambda u^2 dx - \int_{\mathbb{R}^N} F(u)dx,$$

Where $F(u) = \int_0^u f(s)ds$.

1.2. Assumptions and main results

Our approach, aligns with the contemporary discourse in nonlinear analysis, enabling us to probe a broad spectrum of solutions under various nonlinearity conditions. More precisely, for the nonlinearity $f$, we assume that

$(F_1)$ $f \in C^1[0, +\infty), f(0) = 0$ and $f(s) > 0$ for $s > 0$.

$(F_2)$ There exist $\alpha, \beta, \mu_1, \mu_2 > 0$ satisfying

$$2 < \alpha, \beta < 2^* := \frac{2N}{N-2},$$

such that,

$$\lim_{t \to 0^+}\frac{f'(t)}{t^{\alpha-2}} = \mu_1(\alpha - 1) > 0 \text{ and } \lim_{t \to +\infty}\frac{f'(t)}{t^{\beta-2}} = \mu_2(\beta - 1) > 0.$$

$(F_3)$ This assertion articulates the non-existence of positive, radially decreasing classical solutions for $-\Delta u + \frac{k}{2}[\Delta(u^2)]u = f(u)$ in $\mathbb{R}^N$.

REMARK 1.1. By [1] Theorem 2.2-(ii), $(F_3)$ holds for $f(t) = |t|^{p-2}t, p \in (2, 2N/(N-2))$.

Formally this can be related to the question of concerning critical points of a functional $J: Y(c) \to \mathbb{R}$ on the set $Y = \{u \in H^1(\mathbb{R}^N) \mid \int_{R^N}|\nabla u|^2 u^2\, dx < \infty\}$ having $\lambda$ as lagrange multiplier given by

$$J(u) = \frac{1}{2}\int_{\mathbb{R}^N}(1 - \kappa u^2)|\nabla u|^2 dx - \int_{\mathbb{R}^N} F(u)dx.$$

To introduce our result, let us first define $U$ as the unique positive solution of:

$$-\Delta U + U = \mu_1 U^{\alpha-1} \text{ in } \mathbb{R}^N, U(0) = \max_{x \in \mathbb{R}^N} U(x). \tag{1.7}$$

and $V$ the unique positive solution of:

$$-\Delta V + V = \mu_2 V^{\beta-1} \text{ in } \mathbb{R}^N, V(0) = \max_{x \in \mathbb{R}^N} V(x). \tag{1.8}$$

where $U, V$ are non-degenerated (See [33](Proposition 2.1)).

THEOREM 1.2. For $N \geq 3$ and $\lambda > 0$, we have the following conclusions.

(i) Mass Subcritical Case

For $2 < \alpha, \beta < 2 + \frac{4}{N}$, there exists a positive normalized solution of (1.1) for any given positive constant $c$. This solution is denoted as $(\lambda, u_\lambda) \in (0, +\infty) \times H^1_{rad}(\mathbb{R}^N)$.

(ii) Exactly Mass Critical Case

If $\alpha = \beta = 2 + \frac{4}{N}$, denote:

$$c_* := \min\left\{\|U\|_2^2, \left(\sqrt{\frac{1}{6}}\right)^N \|V\|_2^2\right\}, c^* := \max\left\{\|U\|_2^2, \left(\sqrt{\frac{1}{6}}\right)^N \|V\|_2^2\right\},$$

A positive normalized solution $(\lambda, u_\lambda) \in (0, +\infty) \times H^1_{rad}(\mathbb{R}^N)$ to (1.1) is attainable if $c$ is within the interval $(c_*, c^*)$ and no positive normalized solution if $c > 0$ small or large.

(iii) At Most Mass Critical Case

This case has two sub-cases.

(iii-1) If $2 < \alpha < \beta = 2 + \frac{4}{N}$, a positive normalized solution $(\lambda, u_\lambda) \in (0, +\infty) \times H^1_{rad}(\mathbb{R}^N)$ to (1.1) exists for $0 < c < \left(\sqrt{\frac{1}{6}}\right)^N \|V\|_2^2$ and no positive normalized solution if $c > 0$ large.

(iii-2) If $2 < \beta < \alpha = 2 + \frac{4}{N}$, a positive normalized solution $(\lambda, u_\lambda) \in (0, +\infty) \times H^1_{rad}(\mathbb{R}^N)$ to (1.1) exists if $c > \|U\|_2^2$ and no positive normalized solution if $c > 0$ small.

(iv) Mixed Case

This is also divided into two sub-cases.

(iv-1) For $2 < \alpha < 2 + \frac{4}{N} < \beta < 2^*$, there are at least two distinct positive normalized solutions $(\lambda_i, u_{\lambda_i}) \in (0, +\infty) \times H^1_{rad}(\mathbb{R}^N)$ to (1.1) for small values of $c > 0$ and no positive normalized solution for large values of $c > 0$.

(iv-2) For $2 < \beta < 2 + \frac{4}{N} < \alpha < 2^*$, similar outcomes are observed for (1.1) but for large values of $c > 0$ and no positive normalized solution for small values of $c > 0$.

(v) At Least Mass Critical Case

Here, we have,

(v-1) When $2 + \frac{4}{N} = \alpha < \beta < 2^*$, a positive normalized solution $(\lambda, u_\lambda) \in (0, +\infty) \times H^1_{rad}(\mathbb{R}^N)$ to (1.1) exists if $0 < c < \|U\|_2^2$ and no positive normalized solution for large values of $c > 0$.

(v-2) When $2 + \frac{4}{N} = \beta < \alpha < 2^*$, a solution $(\lambda, u_\lambda) \in (0, +\infty) \times H^1_{\text{rad}}(\mathbb{R}^N)$ exists for (1.1) if $c > \left(\sqrt{\frac{1}{6}}\right)^N \|V\|_2^2$ and no positive normalized solution for small values of $c$.

(vi) Mass Supercritical Case

For $2 + \frac{4}{N} < \alpha, \beta < 2^*$, there is a normalized positive solution $(\lambda, u_\lambda) \in (0, +\infty) \times H^1_{\text{rad}}(\mathbb{R}^N)$ for (1.1) for any positive value of $c$.

Moreover, the positive solution $u_\lambda$ satisfies the condition that for some $k_1 > 0$, the solution is bounded within $k \in (0, k_1)$ and $\max_{x \in \mathbb{R}^N} |u_\lambda(x)| \leq \sqrt{\frac{1}{3k}}$.

### 1.3. Strategy of this study

Initially, the study employs a dual approach, as referenced in [39,32,2] to convert the given equation (1.1) into a semilinear elliptic equation. For any fixed positive $\lambda$, it is demonstrated that this equation admits at least one positive and radially symmetric solution, denoted as $v_\lambda$. Drawing on methodologies similar to those in [33], the study uses a blow-up argument combined with a Liouville theorem to analyze the asymptotic behaviors of the positive solutions as $\lambda \to 0^+$ or $\lambda \to +\infty$. The $L^2$-norms of $G^{-1}(v_\lambda)$ are also scrutinized. Then, we aim to establish a continuous branch of positive solutions for all positive values of $\lambda$. Finally, by applying a continuity argument, the study concludes with the proof of Theorem 1.2.

## 2. Functional setting and preliminaries

The natural energy functional $J$ related to the equation (3.1) defined by $J: Y(c) \to \mathbb{R}$.

$$J(u) = \frac{1}{2} \int_{\mathbb{R}^N} (1 - \kappa u^2)|\nabla u|^2 dx - \int_{\mathbb{R}^N} F(u) dx, \tag{2.1}$$

However, a significant challenge in applying the variational method to study (1.1) arises due to the functional $J$ not being well-defined in general. This issue stems from the integral $\int_{\mathbb{R}^N} \kappa u^2 |\nabla u|^2 dx$ not being finite. Additionally, there's a need to ensure the positiveness of the term $1 - \kappa t^2$.

To address these challenges and prove the main result, the study proposes an approach that involves establishing normalized solutions for a modified quasilinear Schrödinger equation, as described by (2.2):

$$-\text{div}(g^2(u)\nabla u) + g(u)g'(u)|\nabla u|^2 + \lambda u = f(u), \ x \in \mathbb{R}^N \tag{2.2}$$

Subject to the constraint $\int_{\mathbb{R}^N} |u|^2 dx = c$, with $g(t) = \sqrt{1 - \kappa t^2}$ for $|t| < \sqrt{1/(3\kappa)}$ and $\kappa > 0$. Clearly, when the function $g(t) = \sqrt{1 - \kappa t^2}$, equation (2.2) turns into (1.6).

Let us consider the function $g:[0,+\infty) \to \mathbb{R}$ given by,

$$g(t) = \begin{cases} \sqrt{1-\kappa t^2} & \text{if } 0 \le t < \sqrt{\dfrac{1}{3\kappa}} \\ \dfrac{1}{3\sqrt{2\kappa}t} + \sqrt{\dfrac{1}{6}} & \text{if } \sqrt{\dfrac{1}{3\kappa}} \le t \end{cases}.$$

REMARK 2.1. Setting $g(t) = g(-t)$ for all $t \le 0$, it follows that $g \in C^1(\mathbb{R}, (\sqrt{1/6}, 1]), g$ is an even function, increases in $(-\infty, 0)$ and decreases in $[0, +\infty)$.

Indeed, (2.2) represents the Euler-Lagrange equation corresponding to the natural energy functional:

$$J_\kappa(u) = \frac{1}{2}\int_{\mathbb{R}^N} g^2(u)|\nabla u|^2 dx + \frac{1}{2}\int_{\mathbb{R}^N} \lambda dx - \int_{\mathbb{R}^N} F(u(x))dx. \qquad (2.3)$$

But due to the presence of $g$ in $J_\kappa(u)$ some additional difficulties arise. In what follows, let us define $G(t) = \int_0^t g(s)ds$.

We observe that the inverse function $G^{-1}(t)$ not only exists but is also an odd function. Additionally, both $G$ and its inverse $G^{-1}$ are continuous and have continuous second derivatives, as denoted by their inclusion in the class $C^2(\mathbb{R})$.

LEMMA 2.2. The functions $g$ and $G^{-1}$ are characterized by following attributes:

$(g_0)\ g \in C^1\left(\mathbb{R}, \left(\sqrt{\dfrac{1}{6}}, 1\right]\right)$, is even, and $g(0) = 1$.

$(g_1)\ \lim_{t \to +\infty} g(t) = \sqrt{\dfrac{1}{6}}$.

$(g_2)\ \lim_{t \to \infty} tg'(t) = 0$.

$(g_3)\ \lim_{t \to 0} \dfrac{G^{-1}(t)}{t} = 1$.

$(g_4)\ \lim_{t \to \infty} \dfrac{G^{-1}(t)}{t} = \sqrt{6}$.

$(g_5)\ t \le G^{-1}(t) \le \sqrt{6}t$, for all $t \ge 0$.

$(g_6)\ -\dfrac{1}{2} \le \dfrac{t}{g(t)}g'(t) \le 0$, for all $t \ge 0$.

$(g_7)\ \dfrac{1}{g(G^{-1}(v))} \le \dfrac{G^{-1}(t)}{t}$ for all $t > 0$.

Proof

$(g_0)$ follows from Remark 2.1.

$(g_1)$ Let us define the function g(t) piecewise, so we'll consider each part separately:

1. For $0 \leq t < \sqrt{\frac{1}{3\kappa}}$, $g(t) = \sqrt{1 - \kappa t^2}$.
2. For $\sqrt{\frac{1}{3\kappa}} \leq t$, $g(t) = \frac{1}{3\sqrt{2\kappa}t} + \sqrt{\frac{1}{6}}$.

So we only need to consider the the second piece for the limit $t \to \infty$.

$$g(t) = \frac{1}{3\sqrt{2\kappa}t} + \sqrt{\frac{1}{6}}.$$

As $t \to \infty$, the term $\frac{1}{3\sqrt{2\kappa}t}$ approaches 0, since the denominator grows without bound. Therefore, the limit of $g(t)$ as $t \to \infty$ is dominated by the constant term $\sqrt{\frac{1}{6}}$.

$$\lim_{t \to \infty} g(t) = \sqrt{\frac{1}{6}} > 0.$$

($g_3$) We are interested in the derivative of the second part, since the limit is as $t \to \infty$.

$$g'(t) = -\frac{1}{3\sqrt{2\kappa}t^2}.$$

And, $\quad \lim_{t \to \infty} t g'(t) = \lim_{t \to \infty} t(-\frac{1}{3\sqrt{2\kappa}t^2}) = \lim_{t \to \infty} (-\frac{1}{3\sqrt{2\kappa}t}) = 0.$

Proof of ($g_4$) − ($g_6$) can be found in [39] (Lemma 3.1) and ($g_7$) can be found in [32] (Lemma 2.1(6))

Now, setting $v = G(u) = \int_0^u g(s)ds$.

By $J_\kappa(u)$ we obtain the functional:

$$I_\kappa(v) = \frac{1}{2}\int_{\mathbb{R}^N} |\nabla v|^2 dx + \frac{1}{2}\int_{\mathbb{R}^N} \lambda |G^{-1}(v)|^2 dx - \int_{\mathbb{R}^N} F(G^{-1}(v)) dx, \tag{2.4}$$

As it is well-defined, due to Lemma 2.2 and the assumptions on the nonlinearity $F(s)$, of class $C^1$ and being an even functional on $H^1(\mathbb{R}^N)$ then the variational functional (2.4) can be transformed into:

$$\bar{I}_\kappa(v) = \frac{1}{2}\int_{\mathbb{R}^N} |\nabla v|^2 dx - \int_{\mathbb{R}^N} F(G^{-1}(v)) dx. \tag{2.5}$$

The constraint now becomes:

$$S(a) = \{v \in H^1(R^N), \| G^{-1}(v) \|_2^2 = c\}.$$

And,

$$\bar{I}'_\kappa(v)\varphi = \int_{\mathbb{R}^N} \left[\nabla v \nabla \varphi + \lambda \frac{G^{-1}(v)}{g(G^{-1}(v))}\varphi - \frac{f(G^{-1}(v))}{g(G^{-1}(v))}\varphi\right] dx. \tag{2.6}$$

If $v$ is identified as a critical point of the functional $\bar{I}_\kappa$, then the function $u = G^{-1}(v)$ is established as a classical solution to (2.2), as detailed by Alves, Wang, and Shen[39]. Consequently, to locate solutions to equation (2.2), it is effective to focus on determining the existence of solutions to the subsequent equation.

$$-\Delta v + \lambda \frac{G^{-1}(v)}{g(G^{-1}(v))} = \frac{f(G^{-1}(v))}{g(G^{-1}(v))}, \quad x \in \mathbb{R}^N \tag{2.7}$$

with the prescribed mass $\int_{\mathbb{R}^N} |G^{-1}(v)|^2 \, dx = c$.

And corresponding energy functional as defined in (2.5).

REMARK 2.3. Confirming the presence of a non-trivial solution $v$ for equation (2.7) implies that $u = G^{-1}(v)$ serves as a non-trivial solution to (1.1). This conclusion holds true provided the condition $\sup_{\mathbb{R}^N} |u| < \left(\frac{1}{3k}\right)^{1/2}$ is satisfied, as mentioned in [43].

The energy functional associated with problem (2.7), which corresponds to the energy functional $J$, is defined as follows:

$$\bar{I}_k(v) = \frac{1}{2} \int_{\mathbb{R}^N} |\nabla v|^2 \, dx - \int_{\mathbb{R}^N} F(G^{-1}(v)) dx.$$

Given that $g$ is a nondecreasing and positive function, and $\bar{I}_k$ is confirmed to be well-defined in the $H^1(\mathbb{R}^N)$ and falls under the $C^1$-class. Consequently, to discover normalized solutions for equation (1.1), it's adequate to focus on the search for normalized solutions to problem (2.7).

In our forthcoming analysis, we plan to adopt and modify some methodologies from [33] to identify normalized solutions for problem (2.7).

Let,

$$h_\lambda(v) = \frac{f(G^{-1}(v))}{g(G^{-1}(v))} - \lambda \frac{G^{-1}(v)}{g(G^{-1}(v))} + \lambda v.$$

Then (2.7) turns into:

$$-\Delta v + \lambda v = h_\lambda(v), \quad x \in \mathbb{R}^N \tag{2.8}$$

LEMMA 2.4. Under the assumptions $(F_1) - (F_3)$ and $(g_0) - (g_1)$, for any $\lambda > 0$, there hold that

(i) For some $1 < p < 2^* - 1$,

$$\limsup_{s \to +\infty} \frac{h_\lambda(s)}{s^p} < \infty.$$

(ii) $h_\lambda(s) = o(s), s \to 0$.

(iii) There exists $T > 0$ such that $\int_0^T h_\lambda(\tau) d\tau > \frac{\lambda}{2} T^2$.

(iv) $\lim_{s \to +\infty} \frac{h_\lambda(s)}{s^{\beta-1}} = \frac{\mu_2}{\left(\sqrt{\frac{1}{6}}\right)^\beta}$.

(v) $h_\lambda(s) \leq s h'_\lambda(s)$ for $s > 0$ small.

Proof. (i) For $p = \beta - 1$ by the assumption $(F_2)$ we can approximate $f(u)$ for large $u$ as
$$f(u) \approx \mu_2 u^{\beta-1}.$$

Substituting this into $h_\lambda(s)$, we have:
$$h_\lambda(s) \approx \frac{\mu_2 (G^{-1}(s))^{\beta-1}}{g(G^{-1}(s))} - \lambda \frac{G^{-1}(s)}{g(G^{-1}(s))} + \lambda s.$$

Using $(g_4)$ and $(g_1)$, for $s \to \infty$:
$$\limsup_{s \to +\infty} \frac{h_\lambda(s)}{s^p} \approx \limsup_{s \to +\infty} \frac{\mu_2 (G^{-1}(s))^{\beta-1}}{s^p} < \infty.$$

Since $\frac{G^{-1}(s)}{s}$ approaches a constant as $s \to \infty$, and $g(G^{-1}(s))$ approaches a positive value $\sqrt{\frac{1}{6}}$.

Hence the limit is finite.

(ii) Since,
$$h_\lambda(s) = \frac{f(G^{-1}(s))}{g(G^{-1}(s))} - \lambda \frac{G^{-1}(s)}{g(G^{-1}(s))} + \lambda s.$$

Then,
$$\lim_{s \to 0} \frac{h_\lambda(s)}{s} = \lim_{s \to 0} \frac{f(G^{-1}(s))}{sg(G^{-1}(s))} - \lambda \lim_{s \to 0} \frac{G^{-1}(s)}{sg(G^{-1}(s))} + \lambda,$$
$$= \lim_{s \to 0} \frac{f(s)}{s} = 0.$$

(iii) Since,
$$\lim_{t \to +\infty} \frac{\int_0^t h_\lambda(\tau) d\tau}{t^2} = \lim_{t \to +\infty} \frac{h_\lambda(t)}{2t},$$
$$= \frac{1}{2} \left[ \lim_{t \to +\infty} \frac{f(G^{-1}(t))}{tg(G^{-1}(t))} - \lambda \lim_{t \to +\infty} \frac{G^{-1}(t)}{tg(G^{-1}(t))} + \lambda \right].$$

Take $s = G^{-1}(t)$, then we get,
$$t = G(s) = \int_0^s g(\tau) d\tau.$$

And,
$$\lambda \lim_{t \to +\infty} \frac{G^{-1}(t)}{tg(G^{-1}(t))} = \lambda \lim_{s \to +\infty} \frac{s}{G(s)g(s)}.$$

Noting that $\lim_{s \to +\infty} g(s) = \sqrt{\frac{1}{6}}$.

$$\lambda \lim_{s \to +\infty} \frac{s}{G(s)g(s)} = \sqrt{6}\lambda \lim_{s \to +\infty} \frac{s}{\int_0^s g(\tau)d\tau},$$
$$= \sqrt{6}\lambda \lim_{s \to +\infty} \frac{1}{g(s)},$$
$$= 6\lambda.$$

Meanwhile, it follows from $(F_2)$ that,
$$\lim_{t \to +\infty} \frac{f(G^{-1}(t))}{tg(G^{-1}(t))} = \lim_{s \to +\infty} \frac{f(s)}{G(s)g(s)} = \sqrt{6}\lim_{s \to +\infty} \frac{f(s)}{G(s)},$$
$$= 6\lim_{s \to +\infty} \frac{f(s)}{s} = +\infty.$$

Thus we obtain that:
$$\lim_{t \to +\infty} \frac{\int_0^t h_\lambda(\tau)d\tau}{t^2} = \infty.$$

which implies that there exists $T > 0$ such that $\int_0^T h_\lambda(\tau)d\tau > \frac{\lambda}{2}T^2$.

(iv) Set $p = \beta - 1$ and notice that $G(s) \geq s$ for any $s > 0$, we have:

$$\lim_{s \to +\infty} \frac{h_\lambda(s)}{s^p} = \lim_{s \to +\infty} \frac{f(G^{-1}(s))}{s^p g(G^{-1}(s))} + \lim_{s \to +\infty} \frac{\lambda s - \lambda \frac{G^{-1}(s)}{g(G^{-1}(s))}}{s^p},$$
$$= \lim_{t \to +\infty} \frac{f(t)}{(G(t))^p g(t)} - \lambda \lim_{t \to +\infty} \frac{t}{(G(t))^p g(t)},$$
$$= \lim_{t \to +\infty} \frac{f(t)}{(G(t))^p g(t)} = \lim_{t \to +\infty} \frac{f(t)}{\left(\frac{t}{\sqrt{6}}\right)^p \sqrt{\frac{1}{6}}},$$
$$= \lim_{t \to +\infty} \frac{f(t)}{t^p} \cdot (\sqrt{6})^{p+1} = (\sqrt{6})^{p+1} \mu_2.$$

(v) Since,
$$h_\lambda(v) = \frac{f(G^{-1}(v))}{g(G^{-1}(v))} - \lambda \frac{G^{-1}(v)}{g(G^{-1}(v))} + \lambda v.$$

From $(g_0)$ and $(g_3)$ we have $g(0) = 1$, $\lim_{t \to 0} \frac{G^{-1}(t)}{t} = 1$ for $\to 0$.

We can assume that $G^{-1}(s) \approx s$.
this implies $f(G^{-1}(s)) \approx f(s)$ and $g(G^{-1}(s)) \approx g(s) = 1$.
$$h_\lambda(s) \approx f(s) - \lambda(s) + \lambda(s).$$

$$h_\lambda(s) \approx f(s), \tag{2.9}$$

And,
$$sh'_\lambda(s) \approx sf'(s), \tag{2.10}$$

From $(F_2)$,
$$f'(t) \approx \mu_1(\alpha - 1)t^{\alpha-2}.$$
$$\int_0^s f'(t)dt \approx \int_0^s \mu_1(\alpha - 1)t^{\alpha-2}dt.$$
$$f(s) - f(0) \approx \frac{\mu_1(\alpha - 1)s^{\alpha-1}}{(\alpha - 1)}$$

Since $f(0) = 0$,
$$f(s) \approx \mu_1 s^{\alpha-1},$$
$$sf'(s) \approx \mu_1(\alpha - 1)\alpha s^{\alpha-1}.$$

Since $\alpha > 2$, from (2.9) and (2.10) we get,
$$f(s) \leq sf'(s).$$

Consequently,
$$h_\lambda(s) \leq sh'_\lambda(s).$$

According to Lemma 2.4, for every fixed value of positive $\lambda$, the function $h_\lambda$ meets the Berestycki-Lions conditions. Drawing from the results in [3](Theorem 1), it follows that for any $\lambda > 0$, (2.8) has a ground state solution within $H^1(\mathbb{R}^N)$. This solution is both positive and radially symmetric. We can represent this as:
$$\mathcal{S} = \left\{(\lambda, v_\lambda) \in (0, +\infty) \times H^1_{\text{rad}}(\mathbb{R}^N) : v_\lambda > 0, (\lambda, v_\lambda) \text{ solves } \right.$$
$$\left. -\Delta v + \lambda \frac{G^{-1}(v)}{g(G^{-1}(v))} - \frac{f(G^{-1}(v))}{g(G^{-1}(v))} = 0 \text{ in } \mathbb{R}^N \right\}.$$

REMARK 2.5. Following [4](Theorem 2), for any $\lambda > 0$, any positive $C^2$ solution $\bar{u}$ of problem (2.8) that approaches zero as $|x| \to \infty$, is radially symmetric around some point $x_0 \in \mathbb{R}^N$. This means $\bar{u}(x)$ can be expressed as $\bar{u}_0(|x - x_0|)$, with $\frac{\partial \bar{u}_0}{\partial r} < 0$ for $r = |x - x_0| > 0$. For simplicity, we assume $x_0 = 0$. The elliptic estimate further implies that both $\bar{u}$ and $|\nabla \bar{u}|$ decrease exponentially at infinity. The function $h_\lambda(t)$ is defined as $h_\lambda(t) = h_1(t) + h_2(t)$, for $t \geq 0$, where:
$$h_1(t) = \frac{f(G^{-1}(t))}{g(G^{-1}(t))}, h_2(t) = -\lambda \frac{G^{-1}(t)}{g(G^{-1}(t))} + \lambda t.$$

Then by $(F_2), h_1(t) = O(|t|^{\alpha-1})$ as $t \to 0$ and by $(g_0)$ and $(g_3)$,

$$\lim_{t\to 0} h_2(t) = \lambda \lim_{t\to 0}\left[t - \frac{G^{-1}(t)}{g(G^{-1}(t))}\right]$$
$$= \lambda \lim_{t\to 0} t - \lambda \lim_{t\to 0}\frac{G^{-1}(t)}{g(G^{-1}(t))}$$
$$= 0.$$

It follows that $h_\lambda(t) = O(|t|^{\min\{\alpha-1,2\}})$ as $t \to 0$. Moreover,

$$h'_2(t) = \lambda g^{-3}(G^{-1}(t))[g^3(G^{-1}(t)) - g(G^{-1}(t)) + G^{-1}(t)g'(G^{-1}(t))] \geq 0, t \geq 0.$$

This leads to the conclusion that $h_2$ is a nondecreasing function for all $t \geq 0$. Regarding the condition on $h_1$, it is specified as follows: there exists some constant $C > 0$ and a number $p > 1$, such that:

$$|h_1(t_1) - h_1(t_2)| \leq C|t_1 - t_2|/|\log\min(t_1,t_2)|^p, \ t_1, t_2 \in \left[0, \max_{x\in\mathbb{R}^N}\bar{u}(x)\right].$$

To establish [4](Theorem 2), its application is limited to the proof of [4](Lemma 6.3). Given that $\bar{u}$ exhibits exponential decay, it can be verified that [4](Lemma 6.3) remains applicable in our scenario, and consequently, 4](Theorem 2)is also valid.

PROPOSITION 2.6. For $N \geq 3$, with $1 < p < 2^* - 1$ and $\mu > 0$, there indeed exists a distinct positive radial solution $U_p^\mu \in H^1(\mathbb{R}^N)$ to:

$$-\Delta u + u = \mu u^p, \quad u \in H^1(\mathbb{R}^N) \tag{2.11}$$

Furthermore, when linearizing equation (2.11) at $U_p^\mu$, the transformation:
$$\varphi \to -\Delta\varphi + \varphi - \mu p (U_p^\mu)^{p-1}\varphi$$
results in a null kernel within $H^1_{rad}(\mathbb{R}^N)$.

Proof.

The findings mentioned are derived from analogous results credited to [7], concerning the unique positive solution $U_p$ in $H^1(\mathbb{R}^N)$ for the equation:

$$-\Delta u + u = u^p. \tag{2.12}$$

It is noted that $U_p$ is a positive solution to (2.12) if and only if $U_p^\mu$, defined as $U_p^\mu := \mu^{\frac{1}{1-p}}U_p$ is a positive solution of (2.11).

DEFINITION 2.7. (see [21](Definition 1)) Consider $U$ as a nonempty open subset of a Banach space $L$ and let $\Phi \in C^1(U,\mathbb{R})$. If $u_0$ is a critical point of $\Phi$ at a certain level $\in \mathbb{R}$, it is classified as being of mountain pass-type (mp-type) under the following condition:

For every open neighborhood $W$ contained in $U$ around $u_0$, the set $W \cap \{u \in L: \Phi(u) < d\}$ is nonempty and lacks path-connectedness.

## 3. Asymptotic behaviors of positive solutions

In this section, our focus is on examining the asymptotic characteristics of $v_\lambda$ when $\lambda = \lambda_n \to 0^+$ or $\lambda = \lambda_n \to +\infty$.

Echoing the approach in [33], we present a corresponding result for the scenario where $\lambda = \lambda_n \to 0^+$.

LEMMA 3.1. Let $\{v_n\}_{n=1}^\infty \subset \mathcal{S}$ with $\lambda = \lambda_n \to 0^+$, the subsequent assertions are valid.

(i) 
$$\limsup_{n \to +\infty} \|v_n\|_\infty < +\infty.$$

(ii)
$$\liminf_{n \to +\infty} \frac{\|v_n\|_\infty^{\alpha-2}}{\lambda_n} > 0.$$

(iii)
$$\limsup_{n \to +\infty} \frac{\|v_n\|_\infty^{\alpha-2}}{\lambda_n} < +\infty.$$

(iv) Set
$$w_n(x) := \lambda_n^{\frac{1}{2-\alpha}} v_n\left(\frac{x}{\sqrt{\lambda_n}}\right), \tag{3.1}$$

Then $w_n$ satisfies,

$$-\Delta w_n + w_n = \frac{h_{\lambda_n}\left(\lambda_n^{\frac{1}{\alpha-2}} w_n\right)}{\lambda_n^{\frac{\alpha-1}{\alpha-2}}} \text{ in } \mathbb{R}^N.$$

And,

$$w_n \to U \text{ in } H^1(\mathbb{R}^N) \text{ and } C_{r,0}(\mathbb{R}^N), \tag{3.2}$$

Where $U \in C_{r,0}(\mathbb{R}^N)$ and is the singular positive solution of equation (1.7).

Proof. The proof follows a method akin to that outlined in [33](Section 4.1)

LEMMA 3.2. Consider the sequence $\{v_n\}_{n=1}^\infty$ within $\mathcal{S}$, where $\lambda = \lambda_n \to +\infty$. Then, after selecting a suitable subsequence,
$$\liminf_{n \to +\infty} \|v_n\|_\infty = +\infty.$$

Proof. By regularity, for any fixed $n$, $v_n \in C^2(\mathbb{R}^N)$ and we may assume that $v_n(0) = \|v_n\|_\infty$. Given this, set:

$$\overline{v_n}(x) := \frac{1}{v_n(0)} v_n\left(\frac{x}{\sqrt{\lambda_n}}\right).$$

We have,

$$1 = \bar{v}_n(0) \leq -\Delta \bar{v}_n(0) + \bar{v}_n(0),$$
$$= \frac{1}{\lambda_n v_n(0)} h_{\lambda_n}(v_n(0)),$$

$$= \frac{f(G^{-1}(v_n(0)))}{\lambda_n v_n(0)g(G^{-1}(v_n(0)))} - \frac{G^{-1}(v_n(0))}{v_n(0)g(G^{-1}(v_n(0)))} + 1,$$

$$\leq \frac{c(|G^{-1}(v_n(0))|^{\alpha-1} + |G^{-1}(v_n(0))|^{\beta-1})}{\lambda_n v_n(0)} + 1 \tag{3.3}$$
$$- \frac{G^{-1}(v_n(0))}{v_n(0)g(G^{-1}(v_n(0)))},$$

Here we used the fact that for some $c > 0, |f(t)| \leq c(|t|^{\alpha-1} + |t|^{\beta-1}), t \geq 0$.

Let $G^{-1}(v_n(0)) = b_n$, then $v_n(0) = G(b_n) = \int_0^{b_n} g(\tau)d\tau$. Proceeding with a contradiction-based approach, let's assume that the sequence $\{b_n\}$ is bounded. Under this assumption, (3.3) implies that,

$$1 \leq \frac{c(b_n^{\alpha-1} + b_n^{\beta-1})}{\lambda_n G(b_n)} + 1 - \frac{b_n}{G(b_n)g(b_n)}$$
$$\leq \frac{c(b_n^{\alpha-2} + b_n^{\beta-2})}{\lambda_n} + 1 - \frac{b_n}{G(b_n)g(b_n)} \tag{3.4}$$

If $b_n \to 0$, then $v_n(0) \to 0$. Recall that $v_n$ satisfies:
$$-\Delta v_n + \lambda_n \frac{G^{-1}(v_n)}{g(G^{-1}(v_n))} = \frac{f(G^{-1}(v_n))}{g(G^{-1}(v_n))} \quad \text{in } \mathbb{R}^N. \tag{3.5}$$

From $(F_2)$, it follows that:
$$-\Delta v_n + \frac{(\lambda_n - 1)G^{-1}(v_n)}{g(G^{-1}(v_n))} \leq 0 \text{ in } \mathbb{R}^N.$$

Multiplying both sides by $v_n$ and then integrating over $\mathbb{R}^N$, we obtain:
$$\int_{\mathbb{R}^N} |\nabla v_n|^2 \, dx + (\lambda_n - 1) \int_{\mathbb{R}^N} \frac{G^{-1}(v_n)v_n}{g(G^{-1}(v_n))} dx \leq 0, \tag{3.6}$$

which is contradiction, due to $G^{-1}(v_n(x))v_n(x) > 0, x \in \mathbb{R}^N$. Therefore, $\{b_n\}$ is positive and remains bounded away from zero. As $n \to \infty$, in (3.21), we arrive at the following limit:
$$\limsup_{n\to\infty} \frac{b_n}{G(b_n)g(b_n)} \leq 0.$$

This results in a contradiction. Consequently, it follows that, up to a subsequence, $v_n(0) \to \infty$.

LEMMA 3.3. Let $\{v_n\}_{n=1}^{\infty} \subset \mathcal{S}$ with $\lambda = \lambda_n \to +\infty$, then:

$$\liminf_{n \to +\infty} \frac{\|v_n\|_{\infty}^{\beta-2}}{\lambda_n} > 0, \qquad (3.7)$$

Proof. Recall that,

$$-\Delta v_n = \frac{f(G^{-1}(v_n))}{g(G^{-1}(v_n))} - \lambda_n \frac{G^{-1}(v_n)}{g(G^{-1}(v_n))} \quad \text{in } \mathbb{R}^N.$$

Assuming, for simplicity, that $v_n(0) = \|v_n\|_{\infty}, \forall n \in \mathbb{N}$, and given that $-\Delta v_n(0) \geq 0$,

$$\frac{f(G^{-1}(v_n(0)))}{g(G^{-1}(v_n(0)))} - \lambda_n \frac{G^{-1}(v_n(0))}{g(G^{-1}(v_n(0)))} \geq 0.$$

Since $v_n(0) = \|v_n\|_{\infty}$,

$$f\left(G^{-1}(\|v_n\|_{\infty})\right) \geq \lambda_n G^{-1}(\|v_n\|_{\infty}).$$

From $(F_2)$ which yields that:

$$\liminf_{n \to +\infty} \frac{\|v_n\|_{\infty}^{\beta-2}}{\lambda_n} > 0.$$

so (3.7) is proved.

LEMMA 3.4. Let $\{v_n\}_{n=1}^{\infty} \subset \mathcal{S}$ with $\lambda = \lambda_n \to +\infty$, then

$$\limsup_{n \to +\infty} \frac{\|v_n\|_{\infty}^{\beta-2}}{\lambda_n} < +\infty, \qquad (3.8)$$

Proof. To establish this, we use a method of contradiction. Let's suppose that

$$\lim_{n \to +\infty} \frac{\|v_n\|_{\infty}^{\beta-2}}{\lambda_n} = +\infty,$$

Set $k = \frac{\beta-2}{2}$ and

$$\tilde{v}_n(x) = \frac{1}{\|v_n\|_{\infty}} v_n\left(\frac{x}{\|v_n\|_{\infty}^{k}}\right).$$

then $\tilde{v}_n(0) = \|\tilde{v}_n\|_{\infty} = 1$ and

$$-\Delta \tilde{v}_n(x) = -\frac{1}{\|v_n\|_{\infty}^{1+2k}} \Delta v_n\left(\frac{x}{\|v_n\|_{\infty}^{k}}\right).$$

Let $x' = \frac{x}{\|v_n\|_{\infty}^{k}}$. Then,

$$-\Delta \tilde{v}_n(x) = -\frac{1}{\|v_n\|_{\infty}^{1+2k}} \Delta v_n(x').$$

And $v_n(x') = \tilde{v}_n(x)\|v_n\|_{\infty}$,

$$-\Delta \tilde{v}_n(x) = \frac{1}{\|v_n\|_\infty^{1+2k}} \left[ \frac{f\left(G^{-1}(\tilde{v}_n(x)\|v_n\|_\infty)\right)}{g\left(G^{-1}(\tilde{v}_n(x)\|v_n\|_\infty)\right)} - \lambda_n \frac{G^{-1}(\tilde{v}_n(x)\|v_n\|_\infty)}{g\left(G^{-1}(\tilde{v}_n(x)\|v_n\|_\infty)\right)} \right],$$

Here $1 + 2k = \beta - 1$,

$$-\Delta \tilde{v}_n(x) = \frac{1}{\|v_n\|_\infty^{\beta-1}} \left[ \frac{f\left(G^{-1}(\tilde{v}_n(x)\|v_n\|_\infty)\right)}{g\left(G^{-1}(\tilde{v}_n(x)\|v_n\|_\infty)\right)} \right.$$
$$\left. - \lambda_n \frac{G^{-1}(\tilde{v}_n(x)\|v_n\|_\infty)}{g\left(G^{-1}(\tilde{v}_n(x)\|v_n\|_\infty)\right)} \right]. \tag{3.9}$$

According to condition $(F_2)$, the right-hand side of equation (3.9) is in $L^\infty(\mathbb{R}^N)$. Selecting an appropriate subsequence, we can assume that $\tilde{v}_n \to \tilde{v}$ in $C^2_{\text{loc}}(\mathbb{R}^N)$. Consequently, from $F_2$ and $g_1$, $\tilde{v}$ becomes a non-negative bounded solution to:

$$-\Delta \tilde{v} = \mu_2 (\sqrt{6})^\beta \tilde{v}^{\beta-1} \text{ in } \mathbb{R}^N.$$

As a result of [33](Theorem 2.5-(ii)), it turns out that $\tilde{v} \equiv 0$, This, however, stands in contradiction with the fact that $\tilde{v}(0) = 1$.

LEMMA 3.5. Let $\{v_n\}_{n=1}^\infty \subset \mathcal{S}$ with $\lambda = \lambda_n \to +\infty$. Define:

$$w_n(x) := \lambda_n^{\frac{1}{2-\beta}} v_n \left( \frac{x}{\sqrt{\lambda_n}} \right), \tag{3.10}$$

Then, up to a subsequence, $w_n \to V^*$ in $C_{r,0}(\mathbb{R}^N)$ and $H^1(\mathbb{R}^N)$ as $n \to +\infty$, where $V^*(x) = (\sqrt{\frac{1}{6}} V(\sqrt{6}x)$ and $V$ is identified as the unique positive solution to equation (1.8).

Proof. Referencing Lemma 3.3 and Lemma 3.4, it follows that:

$$0 < \liminf_{n \to +\infty} \frac{v_n(0)^{\beta-2}}{\lambda_n} \leq \limsup_{n \to +\infty} \frac{v_n(0)^{\beta-2}}{\lambda_n} < +\infty.$$

This implies that $\{w_n\}$ is uniformly bounded in $L^\infty(\mathbb{R}^N)$. Given that $w_n$ satisfies:

$$-\Delta w_n + w_n = \frac{h_{\lambda_n}\left(\lambda_n^{\frac{1}{\beta-2}} w_n\right)}{\lambda_n^{\frac{\beta-1}{\beta-2}}} \text{ in } \mathbb{R}^N, \tag{3.11}$$

That is,

$$-\Delta w_n = \lambda_n^{\frac{1-\beta}{\beta-2}} \frac{f\left(G^{-1}\left(\lambda_n^{\frac{1}{\beta-2}} w_n\right)\right)}{g\left(G^{-1}\left(\lambda_n^{\frac{1}{\beta-2}} w_n\right)\right)} - \lambda_n^{\frac{1}{2-\beta}} \frac{G^{-1}\left(\lambda_n^{\frac{1}{\beta-2}} w_n\right)}{g\left(G^{-1}\left(\lambda_n^{\frac{1}{\beta-2}} w_n\right)\right)}. \tag{3.12}$$

From $(F_2)$, It can be verified that the right-hand side of (3.12) belongs to $L^\infty(\mathbb{R}^N)$. Due to elliptic regularity, and selecting a suitable subsequence, we can assume that $w_n \to V^*$ in $C^2_{loc}(\mathbb{R}^N)$. $V^*$ then fulfills the following

$$-\Delta V^* + 6V^* = \mu_2(\sqrt{6})^\beta (V^*)^{\beta-1} \text{ in } \mathbb{R}^N, V^*(0) = \max_{x \in \mathbb{R}^N} V^*(x).$$

In a manner analogous to the approach used in [33], the desired result can be achieved.

Now, our attention shifts to exploring the asymptotic behavior of $\|G^{-1}(v_n)\|_2$ as $n \to \infty$.

THEOREM 3.6. (i) Let $\{v_n\}_{n=1}^\infty \subset H^1(\mathbb{R}^N)$ be positive solutions to (2.7) with $\lambda = \lambda_n \to 0^+$. Then $\lim_{n\to+\infty} \|v_n\|_\infty = 0$, $\lim_{n\to+\infty} \|\nabla v_n\|_2 = 0$.

And,
$$\lim_{n\to+\infty} \|G^{-1}(v_n)\|_2 = \begin{cases} 0 & \alpha < 2 + \frac{4}{N} \\ \|U\|_2 & \alpha = 2 + \frac{4}{N} \\ +\infty & \alpha > 2 + \frac{4}{N} \end{cases}.$$

(ii) Let $\{v_n\}_{n=1}^\infty \subset H^1(\mathbb{R}^N)$ be positive solutions to (2.7) with $\lambda = \lambda_n \to +\infty$. Then, $\lim_{n\to+\infty} \|v_n\|_\infty = +\infty$, $\lim_{n\to+\infty} \|\nabla v_n\|_2 = +\infty$.

And,
$$\lim_{n\to+\infty} \|G^{-1}(v_n)\|_2 = \begin{cases} +\infty & \beta < 2 + \frac{4}{N} \\ (\sqrt{6})^{-N/2} \|V\|_2 & \beta = 2 + \frac{4}{N} \\ 0 & \beta > 2 + \frac{4}{N} \end{cases}.$$

Proof.

Echoing the methods used in [33], it is possible to demonstrate the asymptotic behavior of $\|v_n\|_\infty$ and $\|\nabla v_n\|_2$. For the subsequent discussion, the focus will be specifically on $\|G^{-1}(v_n)\|_2^2$.

(i) Define $w_n$ as per (3.1), and according to Lemma 3.1, it follows that,

$$\|G^{-1}(v_n)\|_2^2 = \int_{\mathbb{R}^N} |G^{-1}(v_n(x))|^2 \, dx,$$

$$= \lambda_n^{-\frac{N}{2}} \int_{\mathbb{R}^N} \left| G^{-1}\left(v_n\left(\frac{y}{\sqrt{\lambda_n}}\right)\right) \right|^2 \, dy,$$

$$= \lambda_n^{-\frac{N}{2}} \int_{\mathbb{R}^N} \left( \frac{G^{-1}\left(v_n\left(\frac{y}{\sqrt{\lambda_n}}\right)\right)}{v_n\left(\frac{y}{\sqrt{\lambda_n}}\right)} \right)^2 \left(v_n\left(\frac{y}{\sqrt{\lambda_n}}\right)\right)^2 dy$$

$$= \lambda_n^{-\frac{N}{2}} \int_{\mathbb{R}^N} \left( \frac{G^{-1}\left(v_n\left(\frac{y}{\sqrt{\lambda_n}}\right)\right)}{v_n\left(\frac{y}{\sqrt{\lambda_n}}\right)} \right)^2 \left(\frac{w_n(y)}{\lambda_n^{\frac{1}{2-\alpha}}}\right)^2 dy,$$

$$= \lambda_n^{-\frac{N}{2}+\frac{2}{\alpha-2}} \int_{\mathbb{R}^N} \left( \frac{G^{-1}\left(v_n\left(\frac{y}{\sqrt{\lambda_n}}\right)\right)}{v_n\left(\frac{y}{\sqrt{\lambda_n}}\right)} \right)^2 w_n^2(y) dy,$$

Since from $(g_3)$, we know that $\lim_{t \to 0} \frac{G^{-1}(t)}{t} = 1$,

We have,
$$= \lambda_n^{-\frac{N}{2}+\frac{2}{\alpha-2}} (\| U \|_2^2 + o_n(1)).$$

This leads to the attainment of the desired result.

(ii) For the case of $\lambda = \lambda_n \to +\infty$, let $w_n$ be defined as in (3.10). Given that $(g_3)$,
$$\lim_{s \to +\infty} \frac{G^{-1}(s)}{s} = \lim_{t \to +\infty} \frac{t}{G(t)} = \sqrt{6}.$$

by Lemma 3.5,
$$\|G^{-1}(v_n)\|_2^2 = \lambda_n^{-\frac{N}{2}} \int_{\mathbb{R}^N} \left( \frac{G^{-1}\left(v_n\left(\frac{y}{\sqrt{\lambda_n}}\right)\right)}{v_n\left(\frac{y}{\sqrt{\lambda_n}}\right)} \right)^2 \left(v_n\left(\frac{y}{\sqrt{\lambda_n}}\right)\right)^2 dy,$$

$$= \lambda_n^{-\frac{N}{2}+\frac{2}{\beta-2}} \int_{\mathbb{R}^N} \left( \frac{G^{-1}\left(v_n\left(\frac{y}{\sqrt{\lambda_n}}\right)\right)}{v_n\left(\frac{y}{\sqrt{\lambda_n}}\right)} \right)^2 w_n^2(y) dy,$$

Hence from $(g_3)$,

$$= \lambda_n^{-\frac{N}{2}+\frac{2}{\beta-2}} (6\|V^*\|_2^2 + o_n(1)).$$

This leads to the achievement of the desired result.

## 4 Local uniqueness of positive solutions

In this section, we will establish the uniqueness of positive solutions for (2.7) under the condition that $\lambda > 0$ is either sufficiently small or large.

THEOREM 4.1. For $N \geqslant 3$ and under the assumptions that conditions $(g_0) - (g_1)$ and $(F_1) - (F_2)$ are satisfied, (2.7) possesses at most one positive solution in cases where either $\lambda > 0$ is large and condition $(g_2)$ is also met, or $\lambda > 0$ is small and condition $(F_3)$ is concurrently fulfilled.

Proof.

To provide the proof, we will employ the approach detailed in [33].

If $\lambda > 0$ is small and $(F_1) - (F_3)$ are satisfied, let's assume that problem (2.7) has two distinct families of positive solutions, $w_\lambda^{(1)}$ and $w_\lambda^{(2)}$, as $\lambda \to 0^+$. Consider

$$v_\lambda^{(i)}(\cdot) := \lambda^{-\frac{1}{\alpha-2}} w_\lambda^{(i)}(\cdot/\sqrt{\lambda}), \; i = 1,2$$

Then by Lemma 3.1, $v_\lambda^{(1)}, v_\lambda^{(2)} \in H_{rad}^1(\mathbb{R}^N)$ satisfy:

$$-\Delta v + v = \frac{h_\lambda\left(\lambda^{\frac{1}{\alpha-2}} v\right)}{\lambda^{\frac{\alpha-1}{\alpha-2}}} \text{ in } \mathbb{R}^N.$$

and as $\lambda \to 0^+$, $v_\lambda^{(i)} \to U$ both in $C_{r,0}(\mathbb{R}^N)$ and in $H^1(\mathbb{R}^N)$, $i = 1,2$. Setting:

$$\xi_\lambda := \frac{v_\lambda^{(1)} - v_\lambda^{(2)}}{\|v_\lambda^{(1)} - v_\lambda^{(2)}\|_\infty}.$$

From mean value theorem for any $x \in \mathbb{R}^N$, exists some $\theta(x) \in [0,1]$ such that:

$$\frac{h_\lambda\left(\lambda^{\frac{1}{\alpha-2}} v_\lambda^{(1)}\right)}{\lambda^{\frac{\alpha-1}{\alpha-2}}} - \frac{h_\lambda\left(\lambda^{\frac{1}{\alpha-2}} v_\lambda^{(2)}\right)}{\lambda^{\frac{\alpha-1}{\alpha-2}}}.$$

$$-\Delta\xi_\lambda + \xi_\lambda = \lambda^{-1}h'_\lambda\left(\lambda^{\frac{1}{\alpha-2}}\left[\theta(x)v_\lambda^{(1)}(x) + (1-\theta(x))v_\lambda^{(2)}(x)\right]\right)\xi_\lambda \text{ in } \mathbb{R}^N.$$

where for any $t \in \mathbb{R}$. We have:
$$\lambda^{-1}h'_\lambda(t) = \frac{f'(G^{-1}(t))g(G^{-1}(t)) - f(G^{-1}(t))g'(G^{-1}(t))}{\lambda[g(G^{-1}(t))]^3}$$
$$-\frac{g(G^{-1}(t)) - G^{-1}(t)g'(G^{-1}(t))}{[g(G^{-1}(t))]^3} + 1.$$

Since $v_\lambda^{(i)} \to U$ in $C_{r,0}(\mathbb{R}^N)$, one can see that:
$$\lim_{\lambda\to 0^+}\lambda^{\frac{1}{\alpha-2}}\left[\theta(x)v_\lambda^{(1)}(x) + (1-\theta(x))v_\lambda^{(2)}(x)\right] = 0 \text{ uniformly for } x \in \mathbb{R}^N.$$

From $(F_2)$ which implies that for any $x \in \mathbb{R}^N$,
$$\lim_{\lambda\to 0^+}\lambda^{-1}h'_\lambda\left(\lambda^{\frac{1}{\alpha-2}}\left[\theta(x)v_\lambda^{(1)}(x) + (1-\theta(x))v_\lambda^{(2)}(x)\right]\right) = \mu_1(\alpha-1)U^{\alpha-2}.$$

Therefore, after selecting an appropriate subsequence, $\xi_\lambda \to \xi$ in $C^2_{loc}(\mathbb{R}^N)$. This $\xi$ is a radial bounded solution to the following equation:
$$-\Delta\xi + \xi = (\alpha-1)\mu_1 U^{\alpha-2}\xi.$$
Given that $\|\xi\|_\infty = 1$, standard elliptic estimates indicate that $\xi$ is a strong solution. Considering the decay property of $U$ and applying a comparison principle, it follows that $\xi$ decays exponentially to 0 as $|x| \to \infty$. Consequently, $\xi \in C_{r,0}(\mathbb{R}^N) \cap H^1_{rad}(\mathbb{R}^N)$. At this juncture, Proposition 2.6 leads us to a contradiction.

Now focusing on the case where $\lambda > 0$ is large and condition $(g_2)$ is satisfied: Let's assume that problem (2.7) allows for two distinct families of positive solutions, $w_\lambda^{(1)}$ and $w_\lambda^{(2)}$ as $\lambda \to +\infty$. We proceed by considering:
$$v_\lambda^{(i)}(\cdot) := \lambda^{-\frac{1}{\beta-2}}w_\lambda^{(i)}(\cdot/\sqrt{\lambda}), \; i = 1,2$$
$$\xi_\lambda := \frac{v_\lambda^{(1)} - v_\lambda^{(2)}}{\|v_\lambda^{(1)} - v_\lambda^{(2)}\|_\infty}.$$

Similarly as above, there exists some $\theta(x) \in [0,1]$ such that,
$$-\Delta\xi_\lambda + \xi_\lambda = \lambda^{-1}h'\left(\lambda^{\frac{1}{\beta-2}}\left[\theta(x)v_\lambda^{(1)}(x) + (1-\theta(x))v_\lambda^{(2)}(x)\right]\right)\xi_\lambda \text{ in } \mathbb{R}^N.$$

By Lemma 3.5, $v_\lambda^{(i)} \to V$ in $C_{r,0}(\mathbb{R}^N)$ as $\lambda \to +\infty$ for $i = 1,2$. It's easy to understand that for any $x \in \mathbb{R}^N$:
$$\lim_{\lambda\to+\infty}\lambda^{\frac{1}{\beta-2}}\left[\theta(x)v_\lambda^{(1)}(x) + (1-\theta(x))v_\lambda^{(2)}(x)\right] = +\infty.$$

By $(g_2)$, for any $x \in \mathbb{R}^N$,
$$\lim_{\lambda\to+\infty}\left.\frac{g(G^{-1}(t)) - G^{-1}(t)g'(G^{-1}(t))}{[g(G^{-1}(t))]^3}\right|_{t=\lambda^{\frac{1}{\beta-2}}\left[\theta(x)v_\lambda^{(1)}(x)+(1-\theta(x))v_\lambda^{(2)}(x)\right]} = 6.$$

And,

$$\lim_{\lambda\to+\infty}\frac{f'(G^{-1}(t))g(G^{-1}(t))-f(G^{-1}(t))g'(G^{-1}(t))}{\lambda[g(G^{-1}(t))]^3}\bigg|_{t=\lambda^{\frac{1}{\beta-2}}\left[\theta(x)v_\lambda^{(1)}(x)+(1-\theta(x))v_\lambda^{(2)}(x)\right]}$$

$$=\mu_2(\sqrt{6})^\beta V^{\beta-2}.$$

Therefore, by selecting a suitable subsequence, $\xi_\lambda \to \xi$ in $C^2_{loc}(\mathbb{R}^N)$ and $\xi$ is a radial bounded solution of,

$$-\Delta\xi+6\xi=\mu_2(\sqrt{6})^\beta V^{\beta-2}\xi.$$

Similarly, this leads to a contradiction.

## 5. Existence of a curve of positive solution when $\lambda \to 0^+$

PROPOSITION 5.1. Let $N \geq 3$ and under the assumption that conditions $(F_1)$ and $(F_2)$ are met, for any $\lambda > 0$ there exists a positive radial solution $v_\lambda \in H^1(\mathbb{R}^N)$ that solves equation (2.7). Furthermore, like any solution to (2.7), this $v_\lambda$ adheres to the Pohozaev identity.

$$\frac{N-2}{2}\int_{\mathbb{R}^N}|\nabla v|^2\,dx+\frac{N}{2}\lambda\int_{\mathbb{R}^N}|v|^2\,dx=N\int_{\mathbb{R}^N}H_\lambda(v)dx. \tag{5.1}$$

Where $H_\lambda(s) = \int_0^s h_\lambda(t)dt$.

Moreover, defining the functional $I_\lambda: H^1(\mathbb{R}^N) \to \mathbb{R}$. By,

$$I_\lambda(v):=\frac{1}{2}\int_{\mathbb{R}^N}|\nabla v|^2+\frac{1}{2}\int_{\mathbb{R}^N}\lambda|v|^2-\int_{\mathbb{R}^N}H_\lambda(v)dx$$

We have,

$$I_\lambda(v_\lambda)=m_\lambda:=\inf_{\gamma\in\Gamma_\lambda}\max_{t\in[0,1]}I_\lambda(\gamma(t)), \tag{5.2}$$

Where,

$$\Gamma_\lambda:=\left\{\gamma\in C\left([0,1],H^1(\mathbb{R}^N)\right):\gamma(0)=0,I_\lambda(\gamma(1))<0\right\}. \tag{5.3}$$

In particular,

$$m_\lambda=\frac{1}{N}\|\nabla v_\lambda\|_2^2. \tag{5.4}$$

Proof.

In [8], for $N \geq 3$, the existence of a least action solution for equation (2.7) was proven under broad assumptions that are fulfilled by Lemma 2.4, especially when the function $h_\lambda(s)$ is odd. Subsequently, in [20] extended this to $N \geq 2$ by providing a mountain pass characterization for this least action solution. Let's examine if these results are still applicable to our nonlinearity $h_\lambda(s)$, which is not odd. It's clear that the solution $v_\lambda$ obtained by replacing $h_\lambda(s)$ with $\tilde{h}_\lambda(s)$ where

$$\tilde{h}_\lambda(s) = \begin{cases} h_\lambda(s) & \text{if } s \geq 0 \\ -h_\lambda(-s) & \text{if } s \leq 0 \end{cases},$$

is non-negative, it also qualifies as a least action solution for the functional $I_k(v)$. Notably, since $h_\lambda(s) \geq 0$, it follows that,

$$\tilde{H}_\lambda(s) := \int_0^s \tilde{h}_\lambda(t)\, dt \geq H_\lambda(s), \ \forall s \in \mathbb{R}$$

Hence, recalling from [20], for any positive least action solution, there exists a path $\gamma \in \Gamma_\lambda$ such that,

$$\max_{t \in [0,1]} I_\lambda(\gamma(t)) = m_\lambda.$$

$\gamma(t)(x) > 0$ for all $x \in \mathbb{R}^N$ and $t \in (0,1]$. From this, we can deduce that the mountain pass characterization, as detailed in (5.2)-(5.3), remains applicable even if $h_\lambda(s)$ is not odd. Lastly, it's worth noting that (5.4) can be directly inferred by combining the findings from (5.1) and (5.2).

Remembering the definition of a critical point of mountain pass-type as provided in Definition 2.7, it follows that,

LEMMA 5.2. Assume that (F1)-(F2) hold. Any solution $\mathfrak{w} \in H^1(\mathbb{R}^N)$ to (3.12) which satisfies $I_\lambda(\mathfrak{w}) = m_\lambda$ is of mp-type.

Proof.

Suppose that $\mathfrak{w} \in H^1(\mathbb{R}^N)$ as a critical point of $I_\lambda$ with $I_\lambda(\mathfrak{w}) = m_\lambda$, it is necessary to demonstrate that for any open neighborhood $W \subset H^1(\mathbb{R}^N)$ of $\mathfrak{w}$, the set

$$W_\mathfrak{w}^- := W \cap \{v \in H^1(\mathbb{R}^N) : I_\mu(v) < m_\mu\}.$$

is both nonempty and not path-connected.

Since $W$ is open set, it includes a ball $B(\mathfrak{w}, 4r)$ where:

$$B(\mathfrak{w}, 4r) = \{v \in H^1(\mathbb{R}^N) : \|v - \mathfrak{w}\|_{H^1(\mathbb{R}^N)} < 4r\}.$$

Applying [40](Lemma 4.1) with $\delta = 2r$ and an arbitrarily fixed $M > 0$, we can infer the existence of a constant $T > 0$ and a continuous path $\gamma : [0, T] \to H^1(\mathbb{R}^N)$ that fulfills,

(i) $\gamma(0) = 0, I_\lambda(\gamma(T)) < -1, \max_{t \in [0,T]} I_\lambda(\gamma(t)) = I_\lambda(\mathfrak{w})$;

(ii) $\gamma(\tau) = \mathfrak{w}$ for some $\tau \in (0, T)$, and

$$I_\lambda(\gamma(t)) < I_\lambda(\mathfrak{w}),$$

For any $t \in [0, T]$ such that $\|\gamma(t) - \mathfrak{w}\|_{H^1(\mathbb{R}^N)} \geq 2r$.

It is important to note that, after a reparametrization, we can assume $T > 0$ to be 1. We fix $\tau_1 < \tau < \tau_2$ such that,

$$\|\gamma(\tau_1) - \mathfrak{w}\|_{H^1(\mathbb{R}^N)} = \|\gamma(\tau_2) - \mathfrak{w}\|_{H^1(\mathbb{R}^N)} = 3r.$$

The points $\gamma(\tau_1)$ and $\gamma(\tau_2)$ belong to $W_{\mathfrak{w}}^-$. and cannot be connected within the set $W_{\mathfrak{w}}^-$. Indeed, suppose they could be linked by a path $s:(\tau_1,\tau_2) \mapsto W_{\mathfrak{w}}^-$ with $s(\tau_1) = \gamma(\tau_1)$ and $(\tau_2) = \gamma(\tau_2)$. In this case, by considering the path

$$\tilde{\gamma}(t) = \begin{cases} \gamma(t), & \text{for } t \in [0,\tau_1] \\ s(t), & \text{for } t \in [\tau_1,\tau_2), \\ \gamma(t), & \text{for } t \in [\tau_2,T] \end{cases}$$

we would have that $\tilde{\gamma} \in \Gamma_\lambda$ with,

$$\max_{t \in [0,T]} I_\lambda(\tilde{\gamma}(t)) < m_\lambda.$$

In contradiction with the definition of $m_\lambda$.

LEMMA 5.3. Assume that (F1)-(F3) hold. There exists some $\lambda_0 > 0$ small, such that for any $\lambda \in (0,\lambda_0)$, (2.7) has a unique positive solution $v_\lambda \in H^1(\mathbb{R}^N)$. Furthermore, the map $\lambda \mapsto v_\lambda, \lambda \in (0,\lambda_0)$ is continuous. That is, $\{(\lambda,v_\lambda): \lambda \in (0,\lambda_0)\}$ is a curve in $\mathbb{R} \times H^1_{\text{rad}}(\mathbb{R}^N)$.

Proof.

By integrating the insights from Proposition 5.1 and Theorem 4.1, it becomes clear that there exists a unique positive solution $v_\lambda \in H^1(\mathbb{R}^N)$ for (2.7), specifically when $\lambda > 0$ is sufficiently small. This outcome is further supported under conditions (F1)-(F2) and Lemma 3.2 as indicated in [22](Corollary 3.5), and also in reference to [42](Lemma 19). Given the stipulations of Lemma 2.4, it is possible to identify some $s^* > 0$ and a value of $\gamma > 2$ such that:

$$h_\lambda(s)s \geq \gamma H_\lambda(s), \forall s \in [0,s^*]. \tag{5.5}$$

Furthermore, when condition (F3) is also in effect Lemma 3.1 provides us with the following information:

$$\|v_\lambda\|_\infty \to 0 \text{ as } \lambda \to 0^+,$$

and therefore, we can determine a certain $\lambda_0 > 0$ such that,

$$v_\lambda(x) \leq \|v_\lambda\|_\infty \leq s^*, \forall x \in \mathbb{R}^N, \forall \lambda \in (0,\lambda_0]. \tag{5.6}$$

For any $\lambda^* \in (0,\lambda_0)$, we aim to demonstrate the continuity of $v_\lambda$ at $\lambda = \lambda^*$. Specifically, this means showing that for any sequence $\lambda_n \to \lambda^*$, the corresponding solutions converge, i.e., $v_{\lambda_n} \to v_{\lambda^*}$ in $H^1(\mathbb{R}^N)$. without the loss of generality we may assume that,

$$\frac{\lambda^*}{2} \leq \lambda_n \leq \lambda_0, \forall n \in \mathbb{N}$$

Noting that $v_{\lambda_n} \in H^1(\mathbb{R}^N)$ satisfies,

$$\int_{\mathbb{R}^N} |\nabla v_{\lambda_n}|^2 + \lambda_n \int_{\mathbb{R}^N} |v_{\lambda_n}|^2 = \int_{\mathbb{R}^N} H_\lambda(v_{\lambda_n}) v_{\lambda_n} dx. \tag{5.7}$$

we obtain, combining (5.7) with (5.1) and (5.4)-(5.6),

$$(\gamma N - 1)\lambda_n \|v_{\lambda_n}\|_2^2 \leq (1 + |2-N|\gamma)\|\nabla v_{\lambda_n}\|_2^2 \leq (1 + |2-N|\gamma)Nm_{\lambda_0}.$$

Hence,
$$\|v_{\lambda_n}\|_2^2 \leq \frac{2(1+|2-N|\gamma)N}{(\gamma N-1)\lambda^*}m_{\lambda_0} \text{ and } \|\nabla v_{\lambda_n}\|_2^2 = Nm_{\lambda_n} \leq Nm_{\lambda_0}.$$

By $(g_5)$,
$$\|G^{-1}(v_{\lambda_n})\|_2^2 \leq \frac{12(1+|2-N|\gamma)N}{(\gamma N-1)\lambda^*}m_{\lambda_0} \text{ and } \|\nabla v_{\lambda_n}\|_2^2 = Nm_{\lambda_n} \leq Nm_{\lambda_0}.$$

Utilizing the fact that the function $\lambda \mapsto m_\lambda$ is non-decreasing, as supported by the mountain pass characterization detailed in (5.2)-(5.3), we ascertain that the sequence $\{v_{\lambda_n}\}$ is bounded in $H^1(\mathbb{R}^N)$. Moreover, the compact embedding of the radially symmetric subspace $H^1_{\text{rad}}(\mathbb{R}^N)$ into $L^p(\mathbb{R}^N), 2 < p < 2^*$ combined with the governing equation for $v_n$, allows us to demonstrate that $\{v_{\lambda_n}\}$ is compact in $H^1(\mathbb{R}^N)$. Therefore, up to a subsequence, $v_{\lambda_n} \to v$ in $H^1(\mathbb{R}^N)$, where $v \in H^1_{\text{rad}}(\mathbb{R}^N)$ is a positive solution of (2.7) for $\lambda = \lambda^*$ The principle of uniqueness then leads to the conclusion that $v = v_{\lambda^*}$. This implies that the mapping $\lambda \mapsto v_\lambda$ is continuous at $\lambda = \lambda^*$.

THEOREM 5.4. Assume that (F1)-(F2) and $(g_0) - (g_2)$ hold. For any $\lambda > 0$, let $v_\lambda$ denote the unique solution of (2.7) in $H^1(\mathbb{R})$. Then the map $\lambda \mapsto v_\lambda$ is continuous from $(0, +\infty)$ to $H^1(\mathbb{R})$. In particular $\{(\lambda, v_\lambda): \lambda \in (0, +\infty)\}$ is connected.

Proof.

For any given $\lambda_0 > 0$, we aim to establish the continuity of $v_\lambda$ at $\lambda = \lambda_0$. Specifically, this means showing that for any sequence $\lambda_n \to \lambda_0$, the corresponding solutions $v_{\lambda_n} \to v_{\lambda_0}$ in $H^1(\mathbb{R}^N)$. We propose that $\{v_{\lambda_n}\}$ is bounded in $H^1(\mathbb{R}^N)$.

The boundedness of $\{\|\nabla v_{\lambda_n}\|_2^2\}$ follows straight forwardly from (5.4) and the fact that $m \mapsto m_\lambda$ is non-decreasing. To establish the boundedness of $\{\|G^{-1}(v_{\lambda_n})\|_2^2\}$, let's define, for any fixed $n$, a function $l(s) = -\lambda_n s + h_\lambda(s)$. Recalling [8](Theorem 5), we find that,
$$\max_{x \in \mathbb{R}} v_{\lambda_n}(x) = v_{\lambda_n}(0) = s_0,$$
Where $s_0 > 0$ is the unique value for which $L(s) < 0$ for $s \in (0, s_0)$ and $L(s_0) = 0$. Here $L(s) := \int_0^s l(t)dt$.

Since $v_{\lambda_n}(x) \leq s_0, x \in \mathbb{R}$, we have that $L(v_{\lambda_n}(x)) \leq 0, x \in \mathbb{R}$. Hence,
$$\int_\mathbb{R} H_\lambda(v_{\lambda_n})dx \leq \frac{\lambda_n}{2}\|v_{\lambda_n}\|_2^2. \tag{5.8}$$

Recalling that $v_{\lambda_n}$ must satisfy the Pohozaev identity, i.e.,
$$\lambda_n \int_{\mathbb{R}^N} |v_{\lambda_n}|^2 dx = \int_{\mathbb{R}^N} H_\lambda(v_{\lambda_n})dx + \int_{\mathbb{R}^N} |\nabla v_{\lambda_n}|^2 dx, \tag{5.9}$$

We deduce, combining (5.8) and (5.9), that,

$$\|v_{\lambda_n}\|_2^2 \leq \frac{2}{\lambda_n}\|\nabla v_{\lambda_n}\|_2^2.$$

Since by $(g_5)$, we get,

$$\|G^{-1}(v_{\lambda_n})\|_2^2 \leq \frac{12}{\lambda_n}\|\nabla v_{\lambda_n}\|_2^2.$$

Thus, $\left\{\|G^{-1}(v_{\lambda_n})\|_2^2\right\}$ is also bounded. This means the sequence of solutions $\{v_n\} \subset H^1(\mathbb{R}^N)$ is bounded. Keeping in mind that for each $n \in \mathbb{N}$, $v_{\lambda_n}$ is a decreasing function, and employing [41](Proposition 1.7.1), it is inferred that $\{v_n\}$ is compact in $L^p(\mathbb{R}), \forall 2 < p \leq \infty$. Consequently, one can see that $v_n \to v$ and $G^{-1}(v_{\lambda_n}) \to G^{-1}(v)$ in $H^1(\mathbb{R}^N)$, after selecting a suitable subsequence, where $v \in H^1_{rad}(\mathbb{R}^N)$ is a positive solution of equation (2.7) for $\lambda = \lambda_0$. Similar to the argument presented in the proof of Lemma 5.3, the uniqueness of positive solutions allows us to draw a definitive conclusion.

Since similar to [33] we only sketch it.

let,

$$\mathcal{S} = \{(\lambda, v_\lambda) \in (0, +\infty) \times H^1_{rad}(\mathbb{R}^N) :, (\lambda, v_\lambda) \text{ solves } (2.7), v_\lambda > 0\}.$$

Let's define $\tilde{\mathcal{S}} \subset \mathcal{S}$ as the connected component of $\mathcal{S}$ that contains the solutions $(\lambda, v_\lambda)$ for $\in (0, \lambda_0)$. We'll use $P_1: (0, +\infty) \times H^1_{rad}(\mathbb{R}^N) \to (0, +\infty)$ to denote the projection onto the $\lambda$-component. Our objective is to demonstrate that $P_1(\tilde{\mathcal{S}}) = (0, +\infty)$. For a fixed $\lambda > 0$, let's define the norm $\| v \|_\lambda$, which is determined by,

$$\| v \|_\lambda := (\| \nabla v \|_2^2 + \lambda \| v \|_2^2)^{\frac{1}{2}}.$$

is equivalent to the usual norm $\| v \|_{H^1(\mathbb{R}^N)}$. The gradient of $I_k$ with respect to $\langle \cdot, \cdot \rangle_\lambda$ can be computed as,

$$\nabla I_k(v) = v - (-\Delta + \lambda)^{-1} h_\lambda(v) =: v - \mathbb{T}_\lambda(v).$$

Classical degree theory arguments is utlized to analyze the equation.

$$\mathbb{T}_\lambda(v) = v, \lambda \in (0, +\infty), v \in H^1_{rad}(\mathbb{R}^N).$$

From [33] we know that, the operator $\mathbb{T}_\lambda: H^1_{rad}(\mathbb{R}^N) \to H^1_{rad}(\mathbb{R}^N)$ is completely continuous and $v$ is the radial solution to (2.7) iff $v$ is a fix point of $\mathbb{T}_\lambda$ in $H^1_{rad}(\mathbb{R}^N)$. For $\lambda \in (0, \lambda_o)$ from lemma 2.4, part (v), we can get local fix point index

$$\text{ind}(\mathbb{T}_\lambda, v_\lambda) = \deg_{LS}(id - \mathbb{T}_\lambda, N_\varepsilon(v_\lambda), 0) = -1$$

where $\varepsilon > 0$ is a small, $N_\varepsilon$ represents the $\varepsilon$-neighborhood in ($H^1_{rad}(\mathbb{R}^N)$, and "deg" refers to the Leray-Schauder degree. This definition is well-established and applicable when $v_\lambda$ is an isolated fixed point of $\mathbb{T}_\lambda$ in $H^1_{rad}(\mathbb{R}^N)$.

Now, for any fixed $0 < a < b$, define

$$S(a,b) = \{(\lambda, v_\lambda) \in [a,b] \times H^1_{rad}(\mathbb{R}^N): v_\lambda > 0, (\lambda, v_\lambda) \text{ solves}$$
$$-\Delta v + \lambda \frac{G^{-1}(v)}{g(G^{-1}(v))} - \frac{f(G^{-1}(v))}{g(G^{-1}(v))} = \text{ in } \mathbb{R}^N\}.$$

Applying a similar blow-up technique and an ordinary differential equation (ODE) approach, it is easy to demonstrated that the set $S(a,b)$ is compact in both $C_{r,0}(\mathbb{R}^N)$ and $H^1_{rad}(\mathbb{R}^N)$. Utilizing topological degree theory, and given the compactness of $S(a,b)$ along with the fact that the local fixed point index $S(a,b)$ and $\text{ind}(\mathbb{T}_\lambda, v_\lambda) = -1$ for $\lambda \in (0, \lambda_0)$, it follows that $P_1(\tilde{S}) = (0, \infty)$. This implies the projection of the connected component $\tilde{S}$ onto the $\lambda$-component encompasses the entire positive real line.

## 6. Applications to existence, non-existence and multiplicity of Positive Normalized Solutions

PROPOSITION 6.1. There exist a constant $C_1 > 0$ independent of $k$ such that $\|v_\lambda\|_\infty \le C_1$.

Proof.

For each $m \in N$ and $\beta > 1$, let $A_m = \{x \in \mathbb{R}^N: |v_\lambda|^{\beta-1} \le m\}$ and $B_m = \mathbb{R}^N \setminus A_m$. Define:

$$v_m = \begin{cases} v_\lambda |v_\lambda|^{2(\beta-1)}, & \text{in } A_m \\ m^2 v_\lambda, & \text{in } B_m \end{cases}.$$

Note that $v_m \in H^1(\mathbb{R}^N)$, $v_m \le |v_\lambda|^{2\beta-1}$ and

$$\nabla v_m = \begin{cases} (2\beta-1)|v_\lambda|^{2(\beta-1)} \nabla v_\lambda, & \text{in } A_m \\ m^2 \nabla v_\lambda, & \text{in } B_m \end{cases}. \tag{6.1}$$

Using $v_m$ as a test function in (2.7), we deduce that,

$$\int_{\mathbb{R}^N} \left[\nabla v_\lambda \nabla v_m + \lambda \frac{G^{-1}(v_\lambda)}{g(G^{-1}(v_\lambda))} v_m\right] dx = \int_{\mathbb{R}^N} \frac{f(G^{-1}(v_\lambda))}{g(G^{-1}(v_\lambda))} v_m dx, \tag{6.2}$$

By (6.2),

$$\int_{\mathbb{R}^N} \nabla v_\lambda \nabla v_m dx = (2\beta-1) \int_{A_m} |v_\lambda|^{2(\beta-1)} |\nabla v_\lambda|^2 dx + m^2 \int_{B_m} |\nabla v_\lambda|^2 dx, \tag{6.3}$$

Let,

$$w_m = \begin{cases} v_\lambda |v_\lambda|^{\beta-1}, & \text{in } A_m \\ m v_\lambda, & \text{in } B_m \end{cases}.$$

Then $w_m^2 = v_\lambda v_m \le |v_\lambda|^{2\beta}$, and

$$\nabla w_m = \begin{cases} \beta |v_\lambda|^{\beta-1} \nabla v, & \text{in } A_m \\ m \nabla v_\lambda, & \text{in } B_m \end{cases}.$$

Hence,

$$\int_{\mathbb{R}^N} |\nabla w_m|^2 dx = \beta^2 \int_{A_m} |v_\lambda|^{2(\beta-1)} |\nabla v_\lambda|^2 dx + m^2 \int_{B_m} |\nabla v_\lambda|^2 dx. \tag{6.4}$$

Therefore, from (6.3) and (6.4),

$$\int_{\mathbb{R}^N} (|\nabla w_m|^2 - \nabla v_\lambda \nabla v_m) dx = (\beta - 1)^2 \int_{A_m} |v_\lambda|^{2(\beta-1)} |\nabla v_\lambda|^2 dx. \tag{6.5}$$

Combining (6.2), (6.3) and (6.5), since $\beta > 1$, we have,

$$\int_{\mathbb{R}^N} |\nabla w_m|^2 dx \leq \left[\frac{(\beta-1)^2}{2\beta - 1} + 1\right] \int_{\mathbb{R}^N} \nabla v_\lambda \nabla v_m dx,$$

$$\leq \beta^2 \int_{\mathbb{R}^N} \left[\nabla v_\lambda \nabla v_m + \lambda \frac{G^{-1}(v_\lambda)}{g(G^{-1}(v_\lambda))} v_m\right] dx,$$

$$\leq \beta^2 \int_{\mathbb{R}^N} \frac{f(G^{-1}(v_\lambda))}{g(G^{-1}(v_\lambda))} v_m dx.$$

From $(g_7)$,

$$\frac{1}{g(G^{-1}(v_\lambda))} \leq \frac{G^{-1}(v_\lambda)}{v_\lambda} \leq 1.$$

Considering $f(s) = |s|^{p-2} s$, we get:

$$\int_{\mathbb{R}^N} |\nabla w_m|^2 dx \leq \beta^2 \int_{\mathbb{R}^N} \left(G^{-1}(v_\lambda)\right)^{p-1} v_m dx.$$

By Sobolev inequality and $(g_5)$,

$$\left(\int_{A_m} |w_m|^{2^*} dx\right)^{(N-2)/N} \leq S \int_{\mathbb{R}^N} |\nabla w_m|^2 dx.$$

$$\leq S\beta^2 \int_{\mathbb{R}^N} |v_\lambda|^{p-2} w_m^2 dx.$$

By Hölder inequality, we have,

$$\left(\int_{A_m} |w_m|^{2^*} dx\right)^{(N-2)/N} \leq S\beta^2 \| v_\lambda \|_{2^*}^{p-2} \left(\int_{\mathbb{R}^N} |w_m|^{2p_1} dx\right)^{\frac{1}{p_1}}.$$

where $\frac{1}{p_1} + \frac{p-2}{2^*} = 1$. Since $|w_m| \leq |v_\lambda|^\beta$ in $\mathbb{R}^N$ and $|w_m| = |v_\lambda|^\beta$ in $A_m$, we have,

$$\left(\int_{A_m} |v_\lambda|^{\beta 2^*} dx\right)^{\frac{N-2}{N}} \leq S\beta^2 \| v_\lambda \|_{2^*}^{p-2} \left(\int_{\mathbb{R}^N} |v_\lambda|^{2\beta p_1} dx\right)^{\frac{1}{p_1}}.$$

By the Monotone Convergence Theorem, let $m \to \infty$, we have,

$$\| v_\lambda \|_{\beta 2^*} \leq \beta^{\frac{1}{\beta}} \{S \| v_\lambda \|_{2^*}^{p-2}\}^{\frac{1}{2\beta}} \| v_\lambda \|_{2\beta p_1}. \tag{6.6}$$

Setting $\sigma = 2^*/(2p_1)$ and $\beta = \sigma$ in (6.6), we obtain $2p_1 \beta = 2^*$ and

$$\| v_\lambda \|_{\sigma 2^*} \leq \sigma^{\frac{1}{\sigma}} \{S \| v_\lambda \|_{2^*}^{p-2}\}^{\frac{1}{2\sigma}} \| v_\lambda \|_{2^*}, \tag{6.7}$$

Taking $\beta = \sigma^2$ in (6.6), we have,

$$\| v_\lambda \|_{\sigma^2 2^*} \leq \sigma^{\frac{2}{\sigma^2}} \{S \| v_\lambda \|_{2^*}^{p-2}\}^{\frac{1}{(2\sigma^2)}} \| v_\lambda \|_{\sigma 2^*} \tag{6.8}$$

From (6.7) and (6.8), we get:

$$\| v_\lambda \|_{\sigma^2 2^*} \leq \sigma^{\frac{1}{\sigma}+\frac{2}{\sigma^2}} \{S \| v_\lambda \|_{2^*}^{p-2}\}^{\frac{1}{2\left(\frac{1}{\sigma}+\frac{1}{\sigma^2}\right)}} \| v_\lambda \|_{2^*}.$$

Taking $\beta = \sigma^i, i = 1,2,...$, iterating (6.6), we get,

$$\| v_\lambda \|_{\sigma^j 2^*} \leq \sigma^{\sum_{j=1}^{j} \frac{i}{\sigma^i}} \{S \| v_\lambda \|_{2^*}^{p-2}\}^{\frac{1}{2}\sum_{i=1}^{j} \frac{1}{\sigma^i}} \| v_\lambda \|_{2^*}.$$

Therefore, by Sobolev inequality and taking the limit of $j \to +\infty$, we get,

$$\| v_\lambda \|_\infty \leq C_1$$

where $C_1 > 0$ is independent of $k > 0$. This ends the proof.

6.1. Proof of Theorem 1.2

Being similar to [33], Let us introduce the function
$$\tilde{\rho}: \tilde{S} \to \mathbb{R}^+, \ (\lambda, v) \mapsto \|G^{-1}(v)\|_2^2. \tag{6.9}$$

(i) By $P_1(\tilde{S}) = (0, +\infty)$ and Theorem 3.6, if $2 < \alpha < 2 + \frac{4}{N}$, then there exists $(\lambda_n, v_n) \subset \tilde{S}$ with $\lambda_n \to 0^+$ and $\|G^{-1}(v_n)\|_2^2 \to 0$. Similarly, if $2 < \beta < 2 + \frac{4}{N}$, there exists $(\lambda_n', v_n') \subset \tilde{S}$ with $\lambda_n' \to +\infty$ and $\|G^{-1}(v_n')\|_2^2 \to +\infty$. Since $\tilde{S}$ is connected, for any given $c > 0$, there exists $(\lambda_c, \bar{u}_c) \in \tilde{S}$ such that $\tilde{\rho}(\bar{u}_c) = c$, that is, (2.7) possesses a positive normalized solution.

(ii) By
$$\rho(S) \supset \rho(\tilde{S}) \supset (c_*, c^*),$$
we see that for any $c \in (c_*, c^*)$, (2.7) possesses at least one normalized solution $(\lambda, v_\lambda)$ with $\lambda > 0$ and $0 < v_\lambda \in H^1_{rad}(\mathbb{R}^N)$. Recalling Theorem 4.1, (2.7) has a unique solution $v_\lambda > 0$ for $\lambda > 0$ small or large enough. So for any $\delta \in \left(0, \min\left\{\| U \|_2^2, \left(\frac{1}{(6)^{\frac{N}{2}}}\right) \| V \|_2^2\right\}\right)$, there exists some $\Lambda_1 > 0$ small enough and $\Lambda_2 > 0$ large enough such that,

$$\begin{cases} \rho(\lambda, v_\lambda) \in (\| U \|_2^2 - \delta, \| U \|_2^2 + \delta), & \forall \lambda \in (0, \Lambda_1) \\ \rho(\lambda, v_\lambda) \in \left(\left(\sqrt{\frac{1}{6}}\right)^N \| V \|_2^2 - \delta, \left(\sqrt{\frac{1}{6}}\right)^N \| V \|_2^2 + \delta\right), & \forall \lambda \in (\Lambda_2, +\infty) \end{cases}$$

On the other hand, by similar to [33] (Corollary 3.2 and Lemma 3.3), For $0 < \Lambda_1 \leq \Lambda_2, +\infty$, we define the set $V_{\Lambda_2}^{\Lambda_1} \coloneqq \{v \in H^1_{rad}(\mathbb{R}^N): v \text{ is a non negative solution to (2.7) with } \lambda \in [\Lambda_1 \Lambda_2]\}$.

We can find some $M_1, M_2 > 0$ such that,
$$\left\{\| G^{-1}(v) \|_2^2 : v \in V_{\Lambda_1}^{\Lambda_2}\right\} \subset [M_1, M_2]$$

Then for any $c \in (0, c_*) \cup (c^*, +\infty)$, (2.7) has no positive normalized solution.

(iii-1) By
$$\rho(\mathcal{S}) \supset \rho(\tilde{\mathcal{S}}) \supset \left(0, \left(\sqrt{\frac{1}{6}}\right)^N \| V \|_2^2\right),$$

applying a similar argument as (ii), we can prove that (2.7) with constraint has at least one normalized solution $(\lambda, v_\lambda)$ with $\lambda > 0$ and $0 < v_\lambda \in H_{rad}^1(\mathbb{R}^2)$ if $0 < c < \left(\sqrt{\frac{1}{6}}\right)^N \| V \|_2^2$.

And (2.7) do not possess any positive normalized solution provided $c > 0$, large This finishes the proof of (iii-1). The case of (iii-2) will be proved in a similar manner.

(iv-1) In such case, we have $\|G^{-1}(v_\lambda)\|_2^2 \to 0$ both as $\lambda \to 0^+$ and as $\lambda \to +\infty$. Define
$$c_1 := \max_{\tilde{\mathcal{S}}} \rho(\lambda, v_\lambda).$$
Then, there exists some $\lambda^* > 0$ with $\|G^{-1}(v_{\lambda^*})\|_2^2 = c_1$. Then for any $c \in (0, c_1)$, there exist some $\lambda_1 \in (0, \lambda^*)$ and $\lambda_2 \in (\lambda^*, +\infty)$ such that,
$$\rho(\lambda_1, v_{\lambda_1}) = \rho(\lambda_2, v_{\lambda_2}) = c.$$

This means that for any $c \in (0, c_1)$, problem (2.7) admits at least two distinct normalized solutions $(\lambda_i, v_i)$ with $\lambda_i > 0$ and $0 < v_i \in H_{rad}^1(\mathbb{R}^2), i = 1, 2$. Furthermore, by employing a reasoning similar to that used in (ii), we can establish the existence of some $M \geq c_1$ such that (2.7) lacks positive normalized solutions when $c > M$. The case described in (iv-2) can be proven through a comparable methodology. It's important to note that in this case, $\|G^{-1}(v_\lambda)\|_2^2 \to +\infty$ both as $\lambda \to 0^+$ and $\lambda \to +\infty$.

(v-1) By
$$\rho(\mathcal{S}) \supset \rho(\tilde{\mathcal{S}}) \supset (0, \| U \|_2^2),$$

we can prove it in a similar way to (iii-1). The case described in (v-2) can be proven in similar way.

(vi) As in the proof of (i), we have that $\rho(\tilde{\mathcal{S}}) = (0, +\infty)$, and the conclusion follows from Theorem 3.6.

Combining the above arguments and by proposition 6.1, the solution $v_\lambda$ of (2.7) satisfies $\|v_\lambda\|_\infty \leq C_1$.
$$\|u_\lambda\|_\infty = \|G^{-1}(v_\lambda)\|_\infty \leq \sqrt{6}\|v_\lambda\|_\infty \leq \sqrt{6}C_1.$$

Now we only need to show that $\sqrt{6}C_1 < \sqrt{\frac{1}{3k}}$.

Hence by defining $k_1$ as $\frac{1}{18C_1^2}$, we ensure that for any $k$ less that $k_1$, the condition $\sqrt{6}C_1 < \sqrt{\frac{1}{3k}}$ holds. $\|u_\lambda\|_\infty$ will always b less than $\sqrt{\frac{1}{3k}}$.

$$\|u_\lambda\|_\infty = \sqrt{6}C_1 < \sqrt{\frac{1}{3k}} \text{ for all } k \in (0, k_1).$$

*Verifying:* Suppose $\sqrt{6}C_1 < \sqrt{\frac{1}{3k}}$ holds, squaring both sides we get $6C_1^2 < \frac{1}{3k}$. This implies that $k < \frac{1}{18C_1^2}$. Thus this inequality holds only if $k < \frac{1}{18C_1^2}$. let $k_1 = \frac{1}{18C_1^2}$.

$$k < k_1 = \frac{1}{18C_1^2}.$$

This implies that $u_\lambda = G^{-1}(v_\lambda)$ is a positive solution of (3.1).

*Impact of k on solution Bounds: Graphical Evidence*

For $C_1 = 1$, this graph clearly demonstrates that for all values of $k$ less than $k_1$, the blue curve $\sqrt{\frac{1}{3k}}$ remains above the red line $\sqrt{6}C_1$. This visual evidence confirms that the maximum value of the solution $u_\lambda$ is indeed less than $\sqrt{\frac{1}{3k}}$ for all $k$ within this range, thereby validating the conditions stipulated by the theorem.

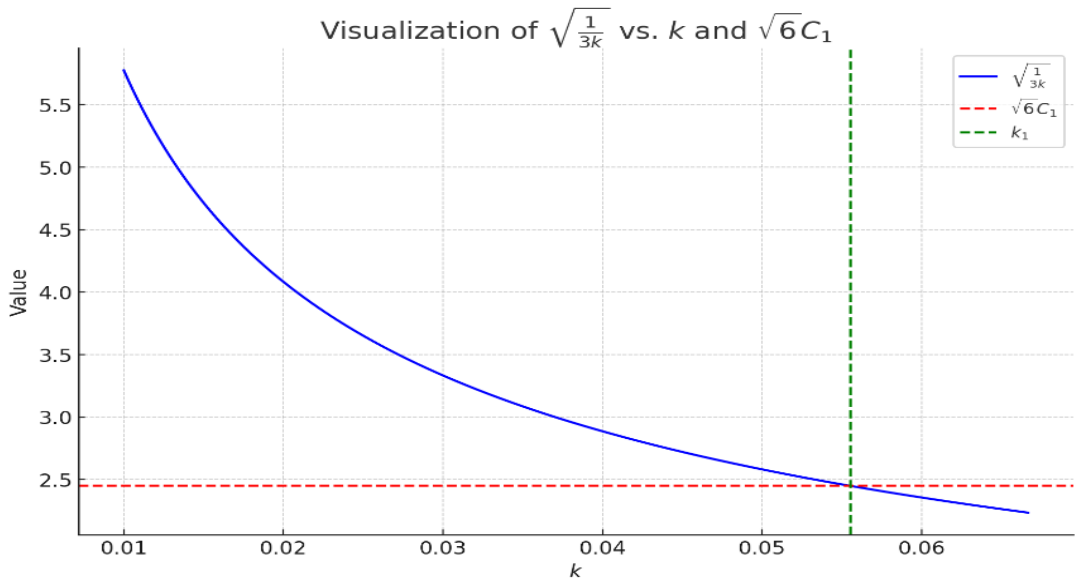

Figure 1 Impact of k on Solution Bound